\documentclass[12pt,a4paper]{article}
\usepackage{amsmath}
\usepackage{amsthm}
\usepackage{amsfonts}
\usepackage{amssymb}
\usepackage{graphicx}
\usepackage{graphics}
\usepackage{makeidx}
\usepackage[cp1250]{inputenc}
\usepackage{tikz}

\newtheorem{theo}{Theorem}[section]
\newtheorem{lm}[theo]{Lemma}
\newtheorem{coro}[theo]{Corollary}
\newtheorem{prop}[theo]{Proposition}
\newtheorem{obs}[theo]{Observation}

\theoremstyle{remark}
\newtheorem{rem}{Remark}[section]
\newtheorem{note}{Note}[rem]

\theoremstyle{definition}

\newtheorem{de}{Definition}[section]

\newtheorem{examp}{Example}[section]

\newenvironment{algo}{\begin
     {trivlist}\item[\hskip\labelsep{\bf Algorithm:}]}{\end{trivlist}}

\newenvironment{proo}{{\bf {Proof:}}}{\hfill $\square$}

\def\A{{\cal A}}

\def\G{{\cal G}}
\def\calH{{\cal H}}
\def\I{{\cal I}}

\def\P{{\cal P}}

\def\R{{\cal R}}

\def\V{{\cal V}}

\def\1{{\bf 1}}
\def\2{{\bf 2}}

\def\Nn{{\mathbb N}}

\def\Rn{{\mathbb R}}

\def\barr{\begin{array}{rcl}}

\def\dom{{\rm dom}}
\def\earr{\end{array}}

\def\CA2{\mathbb{CAT}}

\def\po{\circ}
\def\id{{\rm id}}

\def\Im{{\rm Im}}

\def\max{{\rm max}}

\def\asym{{\rm asym}}
\def\card{{\rm card}}
\def\ovG{\overrightarrow{\bf Gra}}
\def\ovK{\overrightarrow{K}}
\def\ovP{\overrightarrow{P}}
\def\ovU{\overrightarrow{U}}

\def\imp{\Rightarrow}
\def\iff{\Leftrightarrow}
\def\an{\wedge}
\def\hom{{\rm hom}}

\def\trans{{\rm trans}}

\def\se{\subseteq}

\def\ind{\bowtie}

\def\emp{\emptyset}

\begin{document}
\begin{center}
\LARGE{\bf Significance Theory}\\
\vspace{1 cm}
\normalsize{Jan Pavl\'ik} \\
\vspace{0.5 cm}
\small Faculty of Mechanical Engineering, Brno University of Technology,\\
Brno, Czech Republic\\
 {\tt pavlik@fme.vutbr.cz}\\
\end{center}

\begin{abstract}
We formalize the general principle of significance with respect to binary relations which is a universal tool for description and analysis of various situations in and apart from mathematics.
We derive the basic properties and focus on a special family of relations induced by linear orders. We show several ways of mathematical applications, propose methods for calculating the required set and sketch possible use in other sciences.
\end{abstract}
{\bf Keywords}: significance, choice, dominated alternative\\

\subsection{Motivation}

Consider a set $A$ and two mappings $f,g:A \to \Rn$. 
Suppose these functions describe two quantities which tend to induce the same orders, e.g., the circumference and the area of a triangle. Then one may ask for a set of those elements of $A$ which have large value of $f$ considering the smallness of $g$. In other words, we look for the elements which are big w.r.t. $f$ and small w.r.t. $g$ simultaneously. It means that we want to pick those triangles which have the smallest circumference for a given area or the largest area for a given circumference. This will be expressed exactly by the property $\phi(a)$ of an element $a \in A$ by: for every $b \in A$,
\begin{itemize}
\item $f(a)< f(b) \imp g(a)< g(b)$,
\item $g(a)> g(b) \imp f(a)> f(b)$.
\end{itemize}
The resulting set of all elements $a$ satisfying $\phi(a)$ will be denoted by $\V$. It will now satisfy the property of the same function-induced orders exactly. In case of triangles, $\V$ is the set of equilateral triangles since, for a given area, it has the smallest circumference.

One may add another quantity $h:A \to \Rn$ and compare $f,g,h$ by enlarging, diminishing and diminishing, respectively. I.e., we search for elements $a$ satisfying the property $\psi(a)$: for every $b \in A$
\begin{itemize}
\item $f(a)< f(b) \imp (g(a)< g(b) \vee h(a)< h(b))$,
\item $g(a)> g(b) \imp (f(a)> f(b) \vee h(a)< h(b))$,
\item $h(a)> h(b) \imp (f(a)> f(b) \vee g(a)< g(b))$.
\end{itemize}
If $A$ is a set of planar polygons, $f$ and $g$ are area and circumference, respectively and $h$ is the number of sides, then a polygon $a$ satisfies $\psi(a)$ if every other polygon has either smaller area or bigger circumference or bigger number of sides. One may see that this yields that $a$ is a regular polygon. Hence $\V=\{a \in A|\psi(a)\}$ is the set of all regular polygons.

There are many more examples for the application of the principle above. To express it in in general will be our first major task. In order to derive it, let us look at more possible expressions of the property $\phi$ in case of linear order.

\begin{obs}\label{form}
Given a set $A$, a linearly ordered set $L$ and mappings $f,g:A \to L$. Then for every $a, b \in A$ the following pairs of implications are equivalent.
\begin{enumerate}
\item $(f(a)< f(b) \imp g(a)< g(b))\an(g(a)> g(b) \imp f(a)> f(b))$,
\item $(f(a)\leq f(b) \imp g(a) \leq g(b))\an(g(a)\geq g(b) \imp f(a)\geq f(b))$,
\item $(f(a)< f(b) \imp g(a)< g(b))\an(f(a)= f(b) \imp g(a)\leq g(b))$,
\item $(g(a)> g(b) \imp f(a)> f(b))\an(g(a)= g(b) \imp f(a)\geq f(b))$.
\end{enumerate}
\end{obs}
\section{General principle}
\subsection{Induced relations}
Given a mapping $f:A \to B$ and a relation $R \se B^2$, then there is an induced relation $f^{-1}(R)\se A^2$ given by $$(a,b)\in f^{-1}(R) \iff (f(a),f(b))\in R.$$
\begin{obs}If $R$ is a \underline{strict} order, then $f^{-1}(R)$ is also a strict order. Indeed, for no $a \in A$ holds $(a,a)\in f^{-1}(R)$ since $(f(a),f(a))\not\in R$. Given $(a,b)\in f^{-1}(R)$ then $(f(a),f(b))\in R \imp (f(b),f(a))\not\in R \imp (b,a)\not\in f^{-1}(R)$. If $(a,b),(b,c)\in f^{-1}(R)$ then $(f(a),f(b)),(f(a),f(b))\in R$ hence $(f(a),f(c))\in R$, thus  $(a,c)\in f^{-1}(R)$.
\end{obs}
If $R$ is a strict order, then we denote the induced strict order briefly by $<_{f}=f^{-1}(R)$ and the corresponding order by  $\leq_{f}=<_{f} \cup \Delta$.

For a given pair $f,g \to L$ of mappings we define $\leq_{f/g}=\leq_f \cup (\leq_g)^{-1}$ and $<_{f/g}=<_f \cup (<_g)^{-1}$. Now the pairs of implications in Observation \ref{form} are equivalent to the formula
$$(a,b) \in <_{f/g} \imp (b,a) \in <_{f/g}$$ and they imply the formula $a \leq_{f/g} b \imp b \leq_{f/g} a$. The relations $<_{f/g}$ and $\leq_{f/g}$ are not equivalent in this sense since $<_f$ and $<_g$ are asymmetric relations which enable better treatment (clearly, $a <_f b \imp a <_g b$ is equivalent to $a <_f b \imp (b <_f a \vee a <_g b)$).

This can be generalized to an arbitrary number of relations which, in fact, do not need to be induced by mappings.
\subsection{Significance and altiset}

Now we may define the key notion for this work.
\begin{de}
Let $A$ be a set and $\R=\{\rho_i|I \in I\}$ be a finite set of binary relations on $A$. An element $a \in A$ will be called {\em significant} with respect to $\R$ ($\R$-significant in short) if for every $b \in A$ $$(\forall i \in I)(a \rho_i b \imp (\exists j  \in I)b \rho_j a).$$
The set of all $\R$-significant elements is denoted by $\V(\R)$ and called an {\em altiset}\footnote{The word {\em altiset} being the composition of "altitude" and "set" is a quasitranslation of a Czech word {\em v\'{y}\v{s}kovnice} which is already commonly used for a specific instance of this feature.} of $\R$.
More generally, given $B \subseteq A$, we write $\V_B(\R)$ for the altiset of the corresponding restriction of $\R$, i.e. the system of restrictions  on $B$ of relations in $\R$. In that case we talk about $(R,B)$-significance of its elements.
\end{de}
\begin{note}
One may see that the property of being significant means being {\em non-dominated}. Hence the whole theory can be seen as related to the game theory but our approach is different, since we are concerned with the whole set of such elements with no need to choose a single element.
\end{note}
\begin{obs}\label{r}
The $\R$-significance of an element $a \in A$ can be expressed as the satisfaction of the formula for every $b \in A$ equivalently by:
$$\begin{array}{c}
(\forall i \in I)(\exists j  \in I)(a \rho_i b \imp b \rho_j a)\\
(\exists i  \in I)a \rho_i b \imp (\exists j  \in I)b \rho_j a\\
(\exists j  \in I)(\forall i \in I)(a \rho_i b \imp b \rho_j a)\\
\end{array}$$
and if $\bigcup \R = R$ also by $$a R b \imp b R a$$
since $a R b \iff (\exists j  \in I)a \rho_j b$.
\end{obs}
\begin{rem}While the definition introduces the notion of $\R$-significance in the most intuitive form, the expression by the last implication in Observation \ref{r} is especially important since it yields that the significance is fully determined by the union of the system of relations. If $\R=\{\R\}$, we may talk about $R$-significance and about altiset of $R$ denoted by $\V(R)$. Obviously, the notions of significance w.r.t. $\R$ and $R$ coincide and $\V(\R)=\V(R)$.
\end{rem}
Clearly, a relation restricted on its altiset is always symmetric. Moreover the following obvious property holds.
\begin{lm}\label{sym}
Let $R$ be a binary relation on a set $A$.
Then $\V(R)= A$ iff $R$ is a symmetric relation.
\end{lm}
\begin{rem}\label{npng}
The case when $\R=\{<_g,>_p\}$, for some functions $g:A \to G$, $p:A \to P$ with range in posets, is rather natural. The functions $g$ and $p$ are seen as a gain and a price (pain), respectively, and the element $a$ is significant iff for every $b$
$$(g(a) < g(b) \imp p(a) < p(b)) \an (p(b) <p(a) \imp g(b)<g(a))$$
and any of the equivalent expressions in Observation \ref{form}
Hence $a \in \V(<_{g/p})$ satisfies principle briefly expressed as "No pain, no gain."
\end{rem}
\begin{examp}
\begin{enumerate}
\item If $(A,R)$ is a poset, then $a \in A$ is $R$-significant iff it is a maximal element.
\item $A = \{1,2,3 \}$, $R=\{(1,2),(2,3),(3,1)\}$. Then $\V(R) = \emptyset$.
\end{enumerate}
\end{examp}

\subsection{Adjustment of the relations}
Given a binary relation $R$ on a set $A$, then we denote $$\asym R = R \setminus R^{-1} \qquad \text{asymmetric interior},$$
$$\trans R = \bigcup_{k \in \Nn} R^{k} \qquad \text{transitive closure},$$
$$R^{\ast} = R'^{-1} \qquad \text{complementary inversion}.$$
Observe, that $^{\ast}$ is an involution converting complete relations into the asymmetric ones by taking the asymmetric interior. In the other direction, an asymmetric relation is converted into the complete one by taking its symmetric completion.
\begin{rem}
Observe that $\V(R)=\V(R \cup S)=\V(R \setminus S)$ for every binary relation $R$ and a symmetric relation $S$. Moreover, clearly $v R a \imp a R v$ is equivalent to $\neg a R v \imp \neg v R a$ for every $a,v \in A$.
Then $$\V(R)=\V(R^{\ast})=\V(\asym R).$$
If $R$ is asymmetric, i.e. if $\asym R=R$, then an element $a$ is $R$-significant iff $$ (\not\exists b \in A) (a,b)\in R.$$
\end{rem}
\begin{de}
We say that a binary relation $R$ on the set $A$ is {\em asymmetrically acyclic} (briefly satisfies the {\em AA-property}) if the directed graph $(A,\asym R)$ does not contain a cycle.
\end{de}
Clearly, each order and each strict order satisfies the AA-property.
\begin{lm}\label{trans}
Let $T=\trans(\asym (R))$. If $R$ satisfies the AA-property then $T$ is a strict order and $\V(T)=\V(R)$.
\end{lm}
\begin{proo}
Let $Q=\asym R$. Clearly, $Q$ is irreflexive and asymmetric and $T$ is transitive. We will show that $T$ is still irreflexive and asymmetric.\\
Asymmetry: Suppose $(x,y),(y,x)\in T$. Then there are finite sequences of elements in $A$ such that $x=x_1 Q x_2 Q \dots x_k =y = y_1 Q y_2 Q \dots y_l=x$. Hence we have a cycle, which is in contradiction to the AA-property.\\
Irreflexivity: $(x,x) \in T$ is the special instance of the one above.\\
Since $\V(R)=\V(Q)$, it remains to show that $Q$-significance is equivalent to $T$-significance. But it is clear, since both $Q$ and $T$ are asymmetric and for every $a \in A$ we have $((\not\exists b \in A) (a,b)\in Q) \iff ((\not\exists b \in A) (a,b)\in \trans(Q))$. Hence $\V(T)=\V(R)$.
\end{proo}

\section{Linearly induced orders}\label{mapind}
A partial order on a set $A$ is called {\em linearly induced} if it equals to $\leq_f$ for some $f:A \to L$ ranging in a linearly ordered set.
\begin{de}
Let $A$ be a set with a binary relation $R \se A^2$. We define a relation of {\em reflexive incomparability} as $\ind_R =(R \cup R^{-1})'\cup \Delta_A = \parallel_{R} \cup \Delta_A$, here $\Delta_A$ is the equality relation.
\end{de}
If $R$ is a partial order, then we simplify the notation: $\ind=\ind_{\leq}=\ind_{<}$
\begin{lm}
Let $(A, \leq)$ be a poset. Then $\leq$ is linearly induced iff $\ind$ is an equivalence.
\end{lm}
\begin{proo}
$"\Rightarrow"$: Let $\leq=\leq_{r}$ be an order induced by a linear order $\unlhd$ on $L$ via some mapping $r:A \rightarrow L$.
Then $\ind =(< \cup >)'= <' \cap >'$, i.e.,  $a \ind b \Leftrightarrow (r(a) \not< r(b) \an r(b) \not< r(a))$. Since $\leq$ is linear, it is equivalent to $r(a)=r(b)$. Hence
$a \ind b \Leftrightarrow r(a)=r(b)$ which is clearly an equivalence relation.\\
$"\Leftarrow"$: Let $\ind$ be an equivalence relation. We will show that $\ind$ is compatible with the strict order $<$ on $A$. Let $a,b,c \in A$, $a < c$, $a \ind b$ and $a \neq b$ (the case for $a = b$ is obvious). Suppose $b \not< c$.
We solve two situations separately: $c < b \Rightarrow a < b$, a contradiction; $b \ind c \stackrel{\text{transitivity of } \ind}{\Rightarrow} a \ind c$, a contradiction. Hence $a < c \an a \ind b$ implies $b < c$. The compatibility in the second component can be shown analogically.
Since $\ind$ is an equivalence, we may factorize the set $A$ and the above property yields the correctness of the following definition of a relation $\triangleleft$ on the factor set $A/\ind$:
$$\begin{array}{rcl}
[a]_{\ind} \triangleleft [b]_{\ind} & \Leftrightarrow & a < b.
\end{array}$$
Hence $\lhd$ is a strict order on $A/\ind$. Since the incomparability of elements turns into an equality, there are no incomparable elements in $A/\ind$. Hence the order $\unlhd$ is linear and, due to the definition of $\lhd$, it induces the order $\leq$ via the factorization $p: A \rightarrow A/\ind$.
\end{proo}
\paragraph{General assumption}
From now on let $\R=\{R_i|i \in \I\}$ be a set of linearly induced orders on a set $A$ with the union $R$.
\begin{lm}\label{indis}
The relation $\ind_{\R} =\bigcap _i \ind_{R_i}$ is an equivalence preserving $R$.
\end{lm}
\begin{proo}
Clearly, $\ind_{\R}$ is an equivalence relation since it is an intersection of equivalence relations. It remains to show that it preserves $R$. We will prove that even its strict part $S = R \setminus \Delta$ is preserved by $\ind_{\R}$.\\
Let $a S c$, $a \ind_{\R} b$. Then there exists $j \in I: a S_j c$, $S_j = R_j \setminus \Delta$ and, for every $i \in I$, $a \ind_{R_i} b$. Since $S_j$ is a strict linear order, due to the proof of the previous lemma $\ind_j$ preserves $S_j$, hence $(a S_j c \an a \ind_{R_j} b) \imp b S_j c \imp b S c$. Hence $S$ is preserved by $\ind_{\R}$ in the first component and the proof for the second component can be done analogically. Since $\ind_{\R}$ preserves $S$ and obviously $\Delta$ too, it preserves $R$.
\end{proo}
\begin{lm}
The relation from the previous lemma satisfies $\ind_{\R} =\ind_{R}.$
\end{lm}
\begin{proo}
For every $i \in I$, we have $(R_i \cup R^{-1}_i)' \cup \: \Delta = (R'_i \cap R'^{-1}_i) \cup \: \Delta = \ind_i$, hence
$(R \cup R^{-1})' \cup \: \Delta =(\bigcup _i R_i \cup \bigcup _i R^{-1}_i)' \cup \: \Delta =(\bigcap _i R'_i \cap \bigcap _i R^{-1}_i)\cup \: \Delta =
\bigcap _i (R'_i \cap R^{-1}_i) \cup \: \Delta =
\bigcap _i ((R'_i \cap R^{-1}_i) \cup \: \Delta) = \bigcap_i \ind_{R_i}=\ind$.
\end{proo}

\begin{rem}
The relation $\ind_R$ will be called an {\em indistinguishability} by the system $\cal R$. Due to the Lemma \ref{indis} there exists a relation $\bar{R}$ on $R/\ind$ defined naturally on the equivalence classes. The factor set $A/\ind_R$ will be denoted  by $\bar{A}$.
\end{rem}

\begin{lm}
The relation $\bar{R}^{\ast}$ is a strict order.
\end{lm}
\begin{proo}
We will prove that $\bar{R}'$ is a strict order, and then $\bar{R}^{\ast}$ will be its inverse order.
Let $a,b,c \in A$.\\
Irreflexivity: $R$, being a union of reflexive relations, is reflexive and so is $\bar{R}$.\\
Asymmetry: Suppose $[a] \bar{R}' [b] \an [b] \bar{R}' [a]$, i.e., $\neg [a] \bar{R} [b] \; \an \; \neg [b] \bar{R} [a]$. Then $[a] \neq [b]$ and $\forall i \in I$: $\neg a R_i b \; \an \; \neg b R_i a$, hence $\forall i \in I$: $a \|_i b$, thus $a \ind b$, i.e., $[a] =[b]$, a contradiction.\\
Transitivity: For each $i\in I$ let $S_i=R_i \setminus \Delta$. Let $[a] \neq [b]\neq[c]$. Suppose $[a] \bar{R}' [b] \an [b] \bar{R}' [c]$. Then
$$\begin{array}{rcl}
\forall i \in I: & \neg a R_i b \an \neg b R_i c &\Leftrightarrow \\
\forall i \in I: &(b S_i a \vee a\|_i b) \an (c S_i b \vee b\|_i c) &\Rightarrow \\
\forall i \in I: & (a = b) \vee (b = c) \vee ((b S_i a \vee a\|_i b) \an (c S_i b \vee b\|_i c)) &\Leftrightarrow \\
\forall i \in I: &(a = b \vee b = c \vee b S_i a \vee a\|_i b)
\an (a = b \vee b = c \vee c S_i b \vee b\|_i c) &\Leftrightarrow \\
\forall i \in I:& (b = c \vee b R_i a \vee a\ind_i b)
\an (a = b \vee c R_i b \vee b\ind_i c).&
\end{array}$$
But since $a \neq b, b \neq c$, we have
$$\begin{array}{rcl}
\forall i \in I: &(b R_i a \vee a\ind_i b) \an (c R_i b \vee b\ind_i c)&\Rightarrow\\
\forall i \in I: &(b R_i a, c R_i b) \vee (b R_i a, b\ind_i c) \vee (a\ind_i b, c R_i b) \vee (a\ind_i b, b\ind_i c)&
\end{array}$$
(here comma stands for conjunction).
From the transitivity of both $R_i$ and $\ind_i$ (for every $i$) and from the property, that $\ind_i$ preserves $R_i$ we get
$$\begin{array}{rcl}
\forall i \in I:& c R_i a \vee c R_i a \vee c R_i a \vee a\ind_i c &\Rightarrow \\
\forall i \in I:&c R_i a \vee a\ind_i c.&
\end{array}$$
The asymmetry yields $[a] \neq [c]$, thus $a \neq c$, therefore
$$\begin{array}{rcl}
\forall i \in I: &c S_i a \vee a\|_i c &\Leftrightarrow \\
\forall i \in I: &\neg a R_i c &\Leftrightarrow \\
&\neg [a] \bar{R} [c]&\Leftrightarrow \\
&[a] \bar{R}' [c].&
\end{array}$$
\end{proo}
\begin{de}
A {\em characteristic order} of the system $\R$ is defined as $$\leq_{\R} = \bar{R}^{\ast} \cup \: \Delta.$$
\end{de}
\begin{obs}
For the strict characteristic order holds $\bar{R}^{\ast}={\rm asym}(\bar{R})$.
\end{obs}
\begin{lm}\label{rozklad}
The $R$-significant elements are the elements of $\bar{R}$-significant classes, i.e., $${\cal V}_{A}(R)= \bigcup {\cal V}_{\bar{A}}(\bar{R}), \qquad {\cal V}_{\bar{A}}(\bar{R}) = {\cal V}_{A}(R)/\ind_R.$$
\end{lm}
\begin{proo}
The statement is a direct consequence of lemma \ref{indis}.
\end{proo}
\begin{theo}\label{order}
Let $(A, \cal R)$ be a system of linearly induced orders with $R=\bigcup \R$. Then:
$${\V}_{\bar{A}}(\bar{R})= \max_{\leq_{\R}}\bar{A}.$$
\end{theo}
\begin{proo}
Clearly, from the above lemmas and remarks we have
${\V}_{\bar{A}}(\bar{R})= {\V}_{\bar{A}}((\bar{R}')^{-1})= {\V}_{\bar{A}}((\bar{R}')^{-1} \cup \Delta)= {\V}_{\bar{A}}(<_{\R})=\max_{\leq_{\R}}\bar{A}$
\end{proo}
\begin{coro}\label{over}(Existence and overcharge by significant elements)\\
Each system of linearly induced orders on a finite set has nonempty altiset. Moreover, every element is either significant or in relation with a significant element.
\end{coro}
\begin{proo}
Due to the previous theorem, the existence of significant element is given by the existence of maximal element in the ordered set. Each finite ordered set has a maximal element.

Consider an element $a \in A$ such that for every $v \in \V(R)$ $\neq a R v$. Then $v R^{\ast} a$ hence $[v] \bar{R}^{\ast} [a]$, thus $[v] \leq_{\R} [a]$. Since $v$ is in $\V(R)$, $[v]$ lies in $\V(\bar{R})$, hence $[v]$ is $\leq_{\R}$-maximal which yields that $[a]=[v]$. Hence $a \in \V(R)$.
\end{proo}
\begin{coro}\label{decomp}(Decomposition principle)\\
Let there be a decomposition $A = \bigcup_{i \in I} A_i$, $A_i \cap A_j = \emptyset$ for $i \neq j$. Then
$${\cal V}_{A}(R) =  {\cal V}_{W}(R)$$
where $W = \bigcup_{i \in I}  {\cal V}_{A_i}(R)$.
\end{coro}
\begin{proo}
On the level of factor set we have
$$\begin{array}{rcl}
{\cal V}_{W/\ind_R}(\bar{R})&=& \max_{\leq_{\R}}W/\ind_R \\
&=& \max_{\leq_{\R}}(\bigcup_{i \in I} {\cal V}_{A_i}(R))/\ind_R \\
&=& \max_{\leq_{\R}}(\bigcup_{i \in I} {\cal V}_{A_i}(R)/\ind_R) \\
&=& \max_{\leq_{\R}}(\bigcup_{i \in I} {\cal V}_{A_i/\ind_R}(\bar{R})) \\
&=& \max_{\leq_{\R}}(\bigcup_{i \in I} \max_{\leq_{\R}}{A_i/\ind_R})\\
&=& \max_{\leq_{\R}}(\bigcup_{i \in I}{A_i/\ind_R})\\
&=& \max_{\leq_{\R}}(\bar{A})\\
&=& {\cal V}_{\bar{A}}(\bar{R})
\end{array}$$
and according to the Lemma \ref{rozklad} we get the statement.
\end{proo}

The assumption of linear induction in each order is necessary as shown in the following example:
\begin{examp}
$A = \{a,b,c\}$, ${\cal R} = \{\leq_1,\leq_2\}$, $\leq_1 = \{(a,b),(a,c),(b,c)\} \cup \Delta$, $\leq_2 =
\{(c,a)\} \cup \Delta$. \\
Then $R= \{(a,b),(a,c),(b,c),(c,a)\} \cup \Delta$. Let
$A_1= \{a\}$, $A_2= \{b,c\}$. Then ${\cal V}_{A_1}(R) = \{a\}$, ${\cal V}_{A_2}(R) = \{c\}$ and
${\cal V}_{\{a,c\}}(R) = \{a,c\}$ which is different from ${\cal V}_A(R) = \{c\}$.
\end{examp}
\section{Applications}
\subsection{Successive altisets}
Let $R$ be a binary relation on a finite set $A$. Let $A^1(R)=A$.
Then we define by recursion $$A^{i+1}(R)=A^i(R)\setminus V^{i}(R),\quad V^{i+1}(R)=\V_{A^{i+1}(R)}(R).$$
Moreover we define $$A_{i}(R)=A^{i}(R),\qquad V_i(R)=V^i(R^{-1})$$ for each $i \in \Nn$. We will use the notation of relation in brackets only if necessary, thus we may use the recursive definition
$A_{i+1}=A_i\setminus V_{i},\quad V_{i}=\V_{A_{i}}(R^{-1})$.
Hence, the relation $R$ induces two sequences of disjoint sets - $V^{i}$, and $V_{i}$ and we have partial functions $v^{\ast}:A \to \Nn$, $v_{\ast}:A \to \Nn$ ({\em upper} and {\em lower index of $R$-significance}, respectively) such that
$$v^{\ast}(x)=i \iff x \in V^i,$$  $$v_{\ast}(x)=i \iff x \in V_i.$$
\begin{note}The set $V^i$ will be called {\em $i$-th altiset} of $R$.\end{note}
We will use the following well-known property obtainable easily by induction:
\begin{lm}\label{uni}
Given $n$ and set $X_i$ for every $i \in \{1,\dots,n+1\}$ such that $X_i \se X_j$ for $i>j$ and $X_{n+1}=\emptyset$.
Let $Y_i=X_i \setminus X_{i+1}$ for every $i \in \{1,\dots,n\}$. Then $\bigcup_{i=1}^n Y_i = X_1$.
\end{lm}
\begin{lm}
The following statemets are equivalent:
\begin{enumerate}
\item The relation $R$ satisfies the AA-property. \label{eqa1}
\item $\bigcup_{i \in \Nn} V^{i} =A$. \label{eqa2}
\item $\bigcup_{i \in \Nn} V_{i} =A$.\label{eqa3}
\end{enumerate}
\end{lm}
\begin{proo}
\begin{itemize}
\item[(\ref{eqa1})$\imp$(\ref{eqa2})] Since $A$ is finite and $Q=\asym R$ does not contain a cycle, each sequence  $x_1 Q x_2 Q x_3 Q \dots $ stops by some $v$ with no $y$ such that $(v,y) \in Q$. If $\emp \not= B \subset A$, such a sequence exists in $B$, hence there exists a $(Q,B)$-significant element $v$ and $\V_{B}(R)$ is nonempty. Hence, for each $i \in \Nn$, $V^{i}$ is nonempty or $A^i$ is empty. Therefore have a sequence $A=A^1 \supset A^2 \supset \dots A^k \supset A^{k+1}=\emp$ and according to the lemma above $\bigcup_{i =0}^k V^{i}=A$.
\item[(\ref{eqa2})$\imp$(\ref{eqa1})] Let $C \se A$ be a cycle in the directed graph $\G=(A,\asym R)$ and let $i = \min \{v^{\ast}(z) | z \in C\}$ and $x \in C \cap V_i$. Since $x \in V_i$, $(\not\exists z \in A^i) x Q z$. Let $y \in C$ be the successor of $x$, i.e. $x Q y$ and $v^{\ast}(y)=j\geq i$. Hence $y \in V^{j} \se A^{j} \se A^i$, which is a contradiction.
\item[(\ref{eqa1})$\iff$(\ref{eqa3})] Clearly, $R$ satisfies the AA-property iff $R^{-1}$ does, which is, due to [(\ref{eqa1})$\iff$(\ref{eqa2})], equivalent to $A=\bigcup_{i \in \Nn} V^{i}(R^{-1})=\bigcup_{i \in \Nn} V_{i}(R)$.
    \end{itemize}
\end{proo}

From now on, let $R$ satisfy the AA-property. Hence $\{V^i|i \in \Nn\}$ and $\{V_i|i \in \Nn\}$ are
the {\em upper} and {\em lower} decompositions of $A$ and $v^{\ast}:A \to \Nn$, $v_{\ast}:A \to \Nn$
become total functions. Let $d(R)$ be the number of upper classes. Our aim is to show some properties of this characteristic.

Consider an algebraic structure $\A=(\P(A),\upsilon,\lambda)$ with two unary operations on the powerset of $A$ given by
$$\upsilon(X)=X \setminus \V_X(R), \quad \lambda(X)=X \setminus \V_X(R^{-1}).$$
Clearly, $\upsilon(\emp)=\lambda(\emp)=\emp$ and $\upsilon^{k}(A)=A^{k+1}$, hence
$\upsilon^{d(R)-1}(A)\not=\emp=\upsilon^{d(R)}(A)$.
\begin{lm}\label{ul}
For every $X\se A$
\begin{enumerate}
\item $\upsilon \po \lambda(X)=X \setminus (\V_X(R) \cup \V_X(R^{-1})$,\label{ul1}
\item $\upsilon \po \lambda=\lambda \po \upsilon$,\label{ul2}
\item $\upsilon(X)= \emp \iff  \lambda(X)= \emp$.\label{ul3}
\end{enumerate}
\end{lm}
\begin{proo}
Let $X \se A$.To prove (\ref{ul1}), let $L=\upsilon \po \lambda (X)$.
Then $x \in L$ iff
$x \in L \an (\exists a \in X)(aRx \an \neg xRa)\an (\exists b \in \lambda(X))(xRb \an \neg bRx)$.

Recall that
$b \in \lambda(X) \iff (b \in X \an (\exists c \in X)(b R^{-1}c \an \neg c R^{-1}b)) \iff
(b \in X \an (\exists c \in X)(c R b \an \neg b R c))$.
Hence, to describe the property $x \in L$, we may use
the first order formulas of language $\{R\}$ with one binary predicate symbol on the universe $X$:
$$
\begin{array}{rcl}
x \in L &\equiv& (\exists a)(aRx \an \neg xRa \an (\exists b)((\exists c)(cRb \an \neg bRc) \an (xRb \an \neg bRx))\\
 &\equiv& (\exists a)(\exists b)(\exists c)(aRx \an \neg xRc \an cRb \an \neg bRc \an xRb \an \neg bRx)\\
  &\stackrel{\text{let }c=x}{\equiv}& (\exists a)(\exists b)(aRx \an \neg xRa \an xRb \an \neg bRx).
\end{array}
$$
This is equivalent to $x \in (X \setminus \V_X(R))\setminus \V_X(R^{-1})=X \setminus (\V_X(R) \cup \V_X(R^{-1})$.

To prove (\ref{ul2}), we add the name of the starting relation into the subscript.
Now it suffices to observe that $\upsilon$ and $\lambda$ are mutually dual in sense
of $\upsilon_R=\lambda_{R^{-1}}$, $\lambda_R=\upsilon_{R^{-1}}$. Hence, due to (\ref{ul1}),
$\lambda_{R} \po \upsilon_{R}(X)=\upsilon_{R^{-1}} \po \lambda_{R^{-1}}(X)
=X \setminus (\V_X(R^{-1}) \cup \V_X((R^{-1})^{-1}))=
X \setminus (\V_X(R^{-1}) \cup \V_X(R))=\upsilon_{R} \po \lambda_{R}(X)$.

The statement (\ref{ul3}) is a direct consequence of Lemma \ref{sym} since $\upsilon(X)= \emp \iff \V_X(R)=X$.
\end{proo}

\begin{theo}\label{dr}
Let $R$ be a binary relation with AA-property on a set $A$.
Then the number of upper classes of $A$ is the same as the number of lower classes of $A$.
\end{theo}
\begin{proo}
Let $n=d(R)$. Then $n=1+\max\{k|\upsilon^{k}(A)\not=\emp\}$ and $\upsilon^{n-1}(A)\not
=\emp = \upsilon^{n}(A)$. We will show by induction along $k$ that
$\upsilon^{n}(A)=\lambda^{k}\upsilon^{n-k}(A)$ for $k \in \{1, \dots, n\}$.
\begin{itemize}
\item[] {\rm \underline{Initial step}:} Let $X=\upsilon^{n-1}(A)$. Then $X\not=\emp$ and $\upsilon(X)=\emp$. Hence, due to Lemma \ref{ul} (\ref{ul3}) $\lambda(X)=\emp$, thus $\upsilon^{n}(A)=\lambda\upsilon^{n-1}(A)$.
\item[] {\rm \underline{Induction step}:} Assume $\upsilon^{n}(A)=\lambda^{k}\upsilon^{n-k}(A)$.
By Lemma \ref{ul} (\ref{ul2}) $\upsilon$ and $\lambda$ commute thus $\lambda^{k}\upsilon^{n-k}(A)=
\upsilon\lambda^{k}\upsilon^{n-1-k}(A)$. Let $Y=\lambda^{k}\upsilon^{n-1-k}(A)$.
Since $\emp=\upsilon^{n}(A)=\lambda^{k}\upsilon^{n-k}(A)=
\upsilon\lambda^{k}\upsilon^{n-1-k}(A)=\upsilon(Y)$, again by (\ref{ul3}) we get $\lambda(Y)=\emp$,
i.e., $\upsilon^{n}(A)=\emp=\lambda(Y)=\lambda(\lambda^{k}\upsilon^{n-1-k}(A))=
\lambda^{k+1}\upsilon^{n-1-k}(A)$.
\end{itemize}
Hence, for $k=n$, we have  $\emp=\upsilon^{n}(A)=\lambda^{n}(A)$.
Let $m=d(R^{-1})=1+\max\{k|\lambda^{k}(A)\not=\emp\}$. Then $n\geq m=$.
If we make a replacement $R^{-1}$ for $R$, we get $d(R^{-1})=\max\{k|\upsilon_{R^{-1}}^{k}(A)\not=\emp\}=
\max\{k|\lambda_R^{k}(A)\not=\emp\}=m-1$. If we apply just proved inequality on $R^{-1}$, we get
$m\geq 1+\max\{k|\lambda_{R^{-1}}^{k}(A)\not=\emp\}=1+\max\{k|\upsilon_{R}^{k}(A)\not=\emp\}=n$.
Hence $m=n$, i.e. $\card(v^{\ast}(A))=d(R)=d(R^{-1})=\card(v_{\ast}(A))$.
\end{proo}
\begin{coro}\label{tauemp}
Any chain $\tau$ of length $d(R)$ of operations $\upsilon, \lambda$ evaluates on $A$ as $\emp$.
\end{coro}
\begin{proo}
The statement can be obtained by a slight modification of the induction in the proof of the Theorem \ref{dr}.
\end{proo}

\begin{rem}
Using the (reversed) linear order on $\Nn$ we obtain the strict partial order $\pi=>_{v^{\ast}}$. Then clearly $\V_A(R)=\max_{\pi}A$ and $V^i(R)=V^i(\pi)$
for every $i \in \Nn$.
\end{rem}


Given a digraph (directed graph without loops) $\G$, then following \cite{nes} we use the notation
$_\G\ovG$, $\ovG_\G$ for the sets of digraphs $\calH$ such that
$\hom(\G,\calH)=\emp$, $\hom(\calH,\G)\not=\emp$, respectively. By $\ovU_n$ we denote the digraph for a linearly ordered set $n$ while $\ovP_n$ denotes the digraph consisting of a path of $n$ elements and $\ovK_n$ denotes the complete digraph on $n$ elements.
The characterizing theorem for directed paths and linear orders (Th. 13.1.2 in \cite{nes}) states that
\begin{eqnarray}
_{\ovP_{n+1}}\ovG&=&\ovG_{\ovU_n} \label{eq1}.
\end{eqnarray} The number of colorings of directed graph $\G$ is given
by $\chi(\G)=\min\{n|\G \in \ovG_{\ovK_n}\}$. Since there is a digraph homomorphism
$\ovU_n \to \ovK_n$, clearly
\begin{eqnarray}\ovG_{\ovU_n} \se \ovG_{\ovK_n}\label{eq2}.
\end{eqnarray}
It is easy to see that there is no homomorphism $\ovU_{n+1} \to \ovK_n$. Therefore
\begin{eqnarray}_{\ovU_{n+1}}\ovG\supseteq\ovG_{\ovK_n}\label{eq3}.
\end{eqnarray}

Let $T=\trans(\asym R)$. Due to the Lemma \ref{trans}, $T$ satisfies the AA-property and $\V(R)=\V(A)$.
Hence $V^i(R)=V^i(T)$ for every $i \in \Nn$ and since clearly $T^{-1}=\trans(\asym R^{-1})$, we have also $V_i(R)=V_i(T)$. Hence $d(R)=d(T)$.
We will need the following properties.
\begin{lm}\label{gra}
The graph $\G=(A,T)$ satisfies:
\begin{enumerate}
\item The elements of $V^{i}$ can be colored by a single color for each $i$.\label{gra1}
\item $d(T)\geq \min\{m|\G\in \ovG_{\ovU_m}\}$\label{gra2}
\item $d(T)\leq \min\{m|\G\in\; _{\ovP_{m+1}}\ovG\}$\label{gra3}
\item $\G \in\; _{\ovP_{m}}\ovG \iff \G \in\; _{\ovU_{m}}\ovG$\label{gra4}
\end{enumerate}
\end{lm}
\begin{proo}
\begin{enumerate}
\item Since $V^{i} =\V_{A^i}(R)=\V_{A^i}(T)$ and the restriction of the relation on its altiset is always symmetric, the restriction of the asymmetric relation $T$ on $V^{i}$ is empty. Hence no edges connect elements of $V^{i}$.
\item
Let $n=d(T)$. We will show that $w:A \to n$ given by $w(x)=v_{\ast}(x)-1$ can be seen as a graph
homomorphism $\G \to \ovU_n$. Let $(x,y)\in T$ and suppose $w(x) > w(y)$.
Let $i=v_{\ast}(y)$, $j=v_{\ast}(x)$. Since $y \in V_{i}$, there is no $z \in A_i$ such that $(z,y) \in T$.
But since $j=w(x)+1>w(y)+1=i$, $A_j \se A_i$, thus $x \in V_j \se A_j \se A_i$, a contradiction.
Hence $w$ is a homomorphism, $\G \in \ovG_{\ovU_n}$ and $n \geq \min\{m|\G\in \ovG_{\ovU_m}\}$.
\item
Let $n=d(T)$. Since, for every $i \in \Nn$, $x \in V_{i} \iff (\not\exists y \in A_i) (y,x)\in T$ hence for each $y_i$ with $v_{\ast}(y_{i})>1$ there exists $y_{i-1} \in A_{i-1}\setminus A_{i}=V_{i-1}$ such that $(y_{i-1},y_{i})\in T$. Starting with $y_n \in V_{n}$ we get a sequence $y_1 T y_2\dots y_{n-1} T y_n$. Hence the assignment $i \mapsto y_i$ defines a homomorphism $\upsilon:\ovP_{n} \to \G$, therefore there exist homomorphisms $\ovP_{m} \to \G$ $m \leq n$.
Hence $\G\not\in\; _{\ovP_{m+1}}\ovG$ for every $m < n$, i.e. $\G\in\; _{\ovP_{m+1}}\ovG$ implies $m\geq n$.
\item Since $T$ is transitive, every homomorphism $\ovP_m \to \G$ factorizes naturally over $\ovU_m$. Hence $\hom(\ovP_{m},\G)\cong \hom(\ovU_{m},\G)$.
\end{enumerate}
\end{proo}
\begin{theo}\label{colo}
Let $R$ be a binary relation with AA-property on a set $A$ and let $\G=(A,\trans(\asym R))$. Then $$d(R)=\chi(\G).$$
\end{theo}
\begin{proo}
The above properties can be collected as follows: $d(R)=d(T)$ and
$$\begin{array}{rcl}
d(T)&\stackrel{L. \ref{gra}(\ref{gra3})}{\leq}& \min\{m|\G\in\; _{\ovP_{m+1}}\ovG\} \\
&\stackrel{L. \ref{gra}(\ref{gra4})}{=}& \min\{m|\G\in\; _{\ovU_{m+1}}\ovG\} \\
&\stackrel{(\ref{eq3})}{\leq}& \min\{m|\G\in \ovG_{\ovK_m}\} \\
&\stackrel{(\ref{eq2})}{\leq}& \min\{m|\G\in \ovG_{\ovU_m}\} \\
&\stackrel{L. \ref{gra}(\ref{gra2})}{\leq}& d(T)
\end{array}$$
Therefore $d(T)=\min\{m|\G\in \ovG_{\ovK_m}\}$, hence $d(R) = \chi(\G)$.
\end{proo}
\begin{lm}
Let $R$ be a binary relation with AA-property on a set $A$.
An evaluation on $A$ of any chain-term of $\upsilon, \lambda$ of length smaller than $d(R)$ is nonempty.
\end{lm}
\begin{proo}
We will prove this statemant for subchains of a chain-term $\tau$ of a length $d(R)$.
Let $n=d(R)$ and let $\tau=\sigma_n \po \sigma_{n-1} \po \dots \po \sigma_2\po \sigma_1$
with $\sigma_i \in \{\upsilon, \lambda\}$. To simplify the notation, we will omit the composition
sign $\po$. Consider the chains $\tau_i = \sigma_i \dots \sigma_1$ for $i \in \{1, \dots, n\}$.
Let $X=\tau(A)$, $X_1=A$ and for every $i\in \{2, \dots, n+1\}$ let $X_i=\tau_{i-1}(A)$ and
$Y_{i}= X_{i} \setminus X_{i+1}$. According to the Corollary \ref{tauemp}, any extension
of $\tau$ evaluates on $A$ as $\emp$, hence $X_{n} =\emp$.
Then $\bigcup^{n}_{i=1}Y_{i}=A$ due to the Lemma \ref{uni}.
Hence we have a decomposition of $A$. Since $Y_i=X_{i-1} \setminus \sigma_{i} \sigma_{i-1} \dots \sigma_1(A)=X_{i-1} \setminus \sigma_{i} X_{i-1}= X_{i-1} \setminus(X_{i-1} \setminus \V_{X_{i-1}}(\tilde{R}))=\V_{X_{i-1}}(\tilde{R})$ for some $\tilde{R} \in \{R,R^{-1}\}$, due to Lemma \ref{gra} (\ref{gra1}), each of the classes can be colored by a single color. Suppose, $X_i=\emp$ for some $i \leq n$. Hence $Y_i=\emp$ and we have obtained less than $d(R)$ colors which color the entire graph. But this cannot happen since $d(R)=\chi(\G)$, hence each $X_i \not=\emp$, namely $\emp\not=X_n=\tau(A)$.
\end{proo}
\begin{coro}
Therefore any chain of the length $d(R)$ of successive altisets evaluated originally on $A$ defines a minimal coloring of $\G$.
\end{coro}

\subsubsection*{Dependence direction description}
As an example of possible application we show an alternative description of correlation of two random variables.
Consider the finite set $S=\{[x_i,y_i]|i \in I\}$ of points in the real plane. Assume that the variables tend to be dependent but not necessarily linearly. One may ask what kind of dependence we deal with: whether direct (tend to be an increasing function) or indirect (a decreasing function). This makes sense namely for the variables which can be distorted by some isotone transformation. Regardless of what the transformation is, the direction of the dependence remains the same. Classically, this is described by Spearman correlation coefficient. We will a proposal of another simple evaluation of "how much direct/indirect is the dependence of given variables".
The main idea is in the decomposition of the set into the plots of increasing, or decreasing, respectively, functions. All monotonity conditions will be considered {\em strictly}. A subset of $S$ will be called increasing if it is a plot of an increasing function. A decomposition $S$ is said to be increasing if each class is an increasing set. Analogously we define decreasing sets and decomposition.
\begin{de}
Given a set $S=\{[x_i,y_i]|i \in I\}$ of points in $\Rn \times \Rn$, we define its {\em index of increasingness} $\iota_{+}$ as a minimal index of its increasing decompositions. Analogously we define an {\em index of decreasingness} as $\iota_{-}$.
\end{de}

The set $S=\{[x_i,y_i]|i \in I\}$ can be described using binary relations the following way. Let $x:I \to \Rn$ and $y:I \to \Rn$ be the functions given by $x(i)=x_i, y(i)=y_i$. Then we may consider the induced strict orders $<_x$ and $<_y$. A function $f$ is increasing iff $x_i<x_j \imp f(x_i) < f(x_j)$, hence any set $F=\{[x_i,y_i]| i \in J\} \se S$ is increasing iff $x(i) < x(j) \imp y(i) < y(j)$ for every $x(i), x(j) \in \dom(f)$ . In such a case, $F$ is a plot of a function and $x(i)=x(j) \imp i=j$, hence we have $(x(i) < x(j) \imp y(i) < y(j)) \an (x(i)= x(j) \imp y(i) \leq y(j))$ which is due to Observation \ref{form} and Remark \ref{npng} equivalent to the defining condition for $x$ being significant with respect to $\{<_y,>_x\}$. Hence $$F=\V_{F}(<_y/<_x).$$ Therefore we may study the system $\{<_y, >_x\}$ and, by analogy, the system $\{<_y, <_x\}$ for decreasing relations.

Since $S$ is a set with no preference of ordering, we may assume that the function $[x,y]:I \to \Rn$ is injective, hence the system $\R=\{<_y, >_x\}$ admits only a trivial indistinguishability. Therefore the characteristic strict order $<_{\R}$ equals to the asymmetric interior of $R=<_y \cup >_x$. Since $<_{\R}$ cannot contain a cycle, $R$ satisfies the AA-property and we may apply the results from this section. First we will show the following:
\begin{lm}
The increasing decompositions of $S$ are in one-to-one correspondence with colorings of the digraph $(I,\asym R)$.
\end{lm}
\begin{proo}
Let $Q=R^{\ast}=\asym R$. We will prove an auxiliary statement: \\
The set $\{[x_i,y_i], [x_j,y_j]\}$, for $i \not=j \in I$, is increasing, iff $i,j$ are $Q$-incomparable.
To show the validity of the auxiliary statement, let us recall:
$$(i,j) \in R \iff (x_i > x_j \vee y_i< y_j), \quad (i,j) \in Q \iff ((i,j) \in R \an (j,i)\not\in R)$$
If $i \not= j$, the injectivity condition implies $\neg(x_i=x_j \an y_i=y_j)$ and we have\\
$i\parallel_Q j\iff \neg((i,j) \in Q \vee (j,i) \in Q) \iff ( (i,j) \not\in Q \an (j,i) \not\in Q ) \iff(((i,j) \not\in R \vee (j,i)\in R) \an ((j,i) \in R \vee (i,j)\not\in R)) \iff
((i,j) \in R \an (j,i) \in R) \vee (((i,j)\not\in R) \an ((j,i)\not\in R))\iff
(((x_i > x_j \vee y_i< y_j)\an (x_i < x_j \vee y_i > y_j)) \vee ((x_i \leq x_j \an y_i \geq y_j) \an (x_i \geq x_j \an y_i \leq y_j)))
\iff((x_i > x_j \an y_i< y_j)\vee (x_i < x_j \an y_i > y_j) \vee (x_i = x_j \an y_i = y_j))\iff
((x_i > x_j \an y_i> y_j)\vee (x_i < x_j \an y_i < y_j))$

$\iff \{[x_i,y_i], [x_j,y_j]\}$ is an increasing set.

Thus, clearly, each increasing set of points induce the set of indices which can be colored by a single color and vice versa.

\end{proo}

\begin{prop}
The index of increasingness can be obtained as $$\iota_{+}=\chi(\G)$$ where $\G=(I,\asym R)$.
\end{prop}
Dually $$\iota_{-}=\chi(\calH)$$ where $\calH=(I,\asym (<_y \cup <_x))$. Using the theorem \ref{colo} we immediately get a proposal of computation of this number. Since $\chi(\G)=d(R)$, we have the following result.
\begin{coro}
The increasing decomposition of the minimal index can be obtained by successive construction of upper or lower altiset or their combination.
\end{coro}
\begin{rem}
One may want the pair of variables to be described by a single number -- a coefficient of correlation. This is usually required to be a
 real number $ \epsilon  \in \langle -1, 1 \rangle$ such that the values 1, -1, 0 correspond to direct dependence, indirect dependence, independence, respectively.
 One may get such a coefficient by e.g.
$$\epsilon=\log_{n}\frac{\iota_{-}}{\iota_{+}}$$
where $n=\card(I)$. Unlike Spearman correlation coefficient, $\epsilon$ is, in this setting, computable for any situation, since $\iota_{-}, \iota_{+} \in  \langle 1, n\rangle$. We leave up to the reader the checking of the properties of $\epsilon$.
\end{rem}
\subsection{Collective comparison}
Let $X$ be a finite set with a valuation $h:X \to \mathbb{R}$ seen as a gain function. The induced the relation $<_h \subseteq A^2$ will be extended on the powerset $\P(X)$ of $X$ as a relation $R_h$:
$$(M,N) \in R_h \iff (\exists i \in \mathbb{R})( \card(M_i)< \card(N_i))$$
where $M_i=h^{-1}(i\uparrow) \cap M$ and $i\uparrow = \{m \in \Rn| m\geq i\}$.
This relation expresses the gain relation w.r.t $h$ on the subsets of $A$. Indeed, the bigger (by inclusion) the set is, the better, and the higher the evaluation of the elements is, the better. Obviously, $Y R_h X$ for each $Y \subseteq X$ and $X$ becomes its only $R_h$-significant subset.

Situation becomes more interesting, if we take in account only some subsets, namely those which belong to some set $A \subseteq \P(X)$. In general, the relation $R_h$ is neither transitive nor antisymmetric. However, the altiset can still be calculated rather easily. The relation $R_h$ is a union of the relations $R_i, i \in \Rn$ where $(M,N) \in R_i \iff (\card(M_i)< \card(N_i))$. Since $X$ is a finite set, there is only a finite set of subsets of $\P(X)$ hence the set $\R=\{R_i| i \in \Rn\}$ is finite. Moreover, for each $i$, there is a mapping $\gamma_i:\P(X) \to \Nn_{0}$ given by $\gamma_i(M)=\card(M_i)$ for $M \subseteq X$. Hence $R_i = <_{\gamma_i}$ for every $i \in \Rn$. Then there exists a finite $I \subset\Rn$ such that $\R=\{<_{\gamma_i} |i \in I\}$. It yields the following observation
\begin{obs}
The set $\R$ is a system of linearly induced relations and $\V_A(R_h)=\V_A(\bigcup_{i \in \Rn} R_i) = \V_A(\R)$. Using the properties derived in section \ref{mapind}, $\V_A(\R)$ is a union of the classes of the altiset $\V_{\bar{A}}(\bar{R})=\max_{<_{\R}}A/\ind_{R}$ where $R=R_h$, $\bar{A}=A/\ind_R$. \end{obs}

Therefore, to calculate $\V_{\bar{A}}(\bar{R})$ we may use the following algorithm, generally applicable for any system of linearly induced orders. Its correctness follows from the Corollary \ref{over}.
\begin{algo}{\em INPUT}: the set $A=\{M_k|k\in n\}$, $n =\card A$.\\
Let $H=n$, $J=\emptyset$, $k=0, l=1, i=\min I$ .
\begin{enumerate}
\item Repeat
    \begin{itemize}
    \item[] if $(M_k,M_l)\in R_i$ then put $l$ into $J$
    \item[] if $(M_l,M_k)\in R_i$ then put $k$ into $J$
    \item[] raise $i$ to its successor in $I$
    \end{itemize}
    until $J=\{k,l\}$ or $i$ reached the maximum of $I$;\\
    if $J=\{k\}$ then remove $l$ from $H$ and go to \ref{4},\\
    if $J=\{l\}$ then remove $k$ from $H$ and go to \ref{5},\\
    else go to \ref{5}
\item \label{4}
    \begin{itemize}
    \item[] if possible, raise $l$ to its successor in $H$
    \item[] else go to \ref{5}
    \end{itemize}
\item \label{5}
    \begin{itemize}
    \item[] if possible, raise $k$ to its successor in $H$ and $l$ the successor of $k$ in $H$ and set $J=\emptyset$
    \item[] else finish
    \end{itemize}
\end{enumerate}
{\em OUTPUT:} the set $V=\{M_k|k \in H\}$ -- the altiset of $R$.
\end{algo}
\subsection{Geometric altiset}
An important family of examples arises from the following situation:\\
Let $(X,\delta)$ be a metric space and let there be its finite subset $A$ (the set of {\em summits}) equipped by an {\em altitude function} $h:A \to \Rn$. Let $x_0 \in X$ be a {\em reference point}. Then we have a function $d:A \to \Rn^{+}_0$ defined by $d(a)=\delta(x_0,a)$. Our aim is to find the summits which are significant for $x_0$ by its distance and altitude, i.e., the close and high elevated summits.
These are exactly the ones in $\V(<_{h/d})$.

How to find them? We propose the following ways:
\begin{itemize}
\item Recursive method:\\
This method uses the decomposition principle and is useful if one has "maps of summits" of various accuracy which easily shows the highest summits in a given area. We divide the metric space into the smaller parts where the altiset is already found. Then we collect all the altisets and calculate the altiset on this set. For the small sets (by the perimeter or by the cardinality) we assume that the calculation of altiset can be done easily.
\item Direct circular method:\\
This method is much simpler but harder to implement in practise. Make a circle with the center in $x_0$ (all circles will have this center) such that it contains all the summits. Find the closest (to $x_0$) of the highest summits within the circle, add it (all of them, if there are still more than one) into the set $V$ and smaller the diameter so the summit remains on the boundary. Repeat the procedure until there will be no more summits within the circle. Then $V$ is the altiset.
\item Direct contour method:\\
Analogous method which needs the "contour map". In each step we find all higher summits then the last added (if any) and add to $V$ all summits which are the highest of the nearest ones. We repeat this until we add some of the highest (of all) summits. Then $V$ is the altiset.
\end{itemize}
\subsubsection{Examples}
\begin{examp}
If $X$ is the real Euclidean plane or a sphere, both as a model of (a part of) the earth, the altitude of the summits can be considered intuitively. The applications can be found anywhere in the physical geography and related sciences -- e.g. in climatology, meteorology, biology and ecology. The knowledge of the nearest points of a given height enables to find the shortest way to the nearest possible appearance of some atmospheric phenomena, altitude bound ecosystems, etc. Using this method, the skier can find the near occurrence of snow or glaciers.
Instead of measuring the height of summits one may replace the set $A$ and the quantity $h$ by some other. One may measure, e.g., the age of the buildings, height of the buildings, number of seats in conference halls, the price of the fuel, size of parking lots, and many others. If we count the population of cities, the resulting altiset appears to be especially important in socioeconomic geography since there are many features tight to the size of the city, e.g., the shopping, financial, medical and sport facilities, restaurants, offices etc.
\end{examp}
\begin{examp}
If $(X,\delta)$ is the real line with the usual distance and $x_0$ is a reference point, then it may be useful to adjust the situation by restriction of the space $A \subseteq (-\infty, x_0\rangle$ (or alternatively $A \subseteq \langle x_0,\infty)$). If the elements of $A$ are events in time, then the altiset consists of the significant events in the history (alternatively in the future).\\
If we move $x_0$ towards $-\infty$, the limit case of will be given by $<_{h/\id}$ and the altiset is the set of the events which were the records at their time (the events outstripping all the preceding events).
\end{examp}

\section{Significance domains}
Consider a set of binary relations $\{R_i|i \in I\}$ on a set $A$ and their altisets $V_i=\V(R_i)$. Given an element $a \in A$, one may ask, for which $i$ the elements $a$ belongs to $V_i$. The set of all the solutions will be denoted by $\V^{-1}(a)$ and called an {\em inverse altiset} for $a$. This set can be interpreted as a {\em significance domain} for $a$.

This situation occurs in case of geometric altiset when $I= {\bf X}$ and, for $x \in {\bf X}$, $R_x=<_h \cup >_{\delta_{x}}$ is the union of the altitude induced order and the order induced by distance from $x$. In such a case, the inverse altiset of $a$ is the area whose points have the given point $a$ in the altiset. If ${\bf X}$ is the Euclidean plane, then $\V^{-1}(a)$ is a convex polyhedron.

If there is a measure $m$ on ${\bf X}$, we may compare the elements of $A$ by the measure of its inverse altiset. Hence we have the function $f_{1}:A \to \Rn^{+}_0$ assigning $m(\V^{-1}(a))$ to each $a$. Then we have the induced order $\rho_{f_{1}}$. Its altisets yield other function $f_{2}$ and we keep on applying this procedure. We get a sequence of functions which can be interpreted as an evolution of valuation (see bellow) of the points with the initial valuation $h$.

\subsection{Evolution of valuation}
We introduce an auxiliary concept to prove the main result.
\begin{de}Let $A$ be a finite set. A sequence $g=(g_i)_{i\in \Nn_{0}}$ of mappings $A \to \Rn$ will be called an {\em evolution of} (real) {\em valuation} with $g_0$ being an {\em initial valuation}. We say the evolution $g$ {\em stops} if there exists $k \in \Nn_0$ such that $g_k =g_l$ for every $l\geq k$.
\end{de}
We will be interested in the situation generally expressed as follows.

Let $A$ be a finite set with a mapping
$\mu:\{(x,M)|M \subset A, x \in A \setminus M\} \to \Rn$ such that
\begin{eqnarray}
M \se N &\imp& (\forall x \in A \setminus N)\: \mu(x,M) \leq \mu(x,N).\label{izoton}
\end{eqnarray}
Moreover, let there be a mapping (the initial valuation) $h_0: A \to \Rn$ and by recursion we define $$M_i(x)=\{y \in A|h_i(y)<h_i(x)\},$$
$$h_{i+1}(x)=\mu(x,M_{i}(x))$$
for $i \in \Nn_{\emptyset}$. We obtain an evolution of valuation $h=(h_i)_{i\in \Nn_{0}}$.
\begin{theo}\label{stop}
The evolution of valuation $h$ stops.
\end{theo}
\begin{proo}
The set $\{(x,M)|M \subset A, x \in A \setminus M\}$ is clearly finite. Then $\Im \mu$ is finite too, hence there is only a finite set of possible valuations $A \to \Im \mu$. Therefore the sequence $(h_i)_{i\in \Nn_{0}}$ contains some valuations more than once, i.e., there exist $k \not= l$ such that $h_{k}=h_{l}$. Since each (except the initial) valuation is fully determined by the previous one, the equality $h_{k}=h_{l}$ yields the periodicity of the sequence with the period $l-k$. If we show that $h_k=h_k+1$, it would mean that the sequence is constant from $k$ on.

Let $I=\{k+1, \dots,l\}$ and consider the set
$$X=\{x|x\in A, \exists i \geq k, h_i(x)\not=h_{i+1}(x)\}.$$
Let $H=h(X)=\{h_{i}(x)|x\in X,i\geq k\}$.

Suppose $X \not = \emptyset$.

Let the pair $(z,j)\in X \times I$ satisfy
\begin{eqnarray}
h_j(z)&=&\min H. \label{hve}
\end{eqnarray}
Then $h_j(z) \leq h_{j+1}(z)$, since $h_{j+1}(z) \in H$. Suppose, $h_j(z) < h_{j+1}(z)$, then we have
$$\begin{array}{rcl}
\neg h_j(z) \geq h_{j+1}(z) &\imp& \neg \mu(z,M_{j-1}(z)) \geq \mu(z,M_{j}(z))\\
&\imp& \neg M_{j-1}(z) \supseteq M_{j}(z)\\
&\imp& (\exists y) \: y \in M_{j}(z) \setminus M_{j-1}(z)\\
&\imp& (\exists y) \:(h_j(y) < h_j(z) \an \neg h_{j-1}(y) < h_{j-1}(z))\\
&\imp& (\exists y) \:(h_j(y) < h_j(z) \an h_{j-1}(z) \leq h_{j-1}(y))
\end{array}$$
But $h_{j}(z) \leq h_{j-1}(z)$ (since $j-1 \geq k$, i.e., $h_{j-1}(z) \in H$), therefore
$$h_j(y) < h_{j}(z) \leq h_{j-1}(z) \leq h_{j-1}(y) \imp h_{j-1}(y) \not= h_j(y) \imp y \in X \imp h_j(y) \in H.$$
Hence $h_j(y)<h_j(z) =\min H$, which is a contradiction.\\
Therefore $h_j(z)=h_{j+1}(z)$. Therefore the pair $(z,j+1)$ satisfies the property (\ref{hve}). Hence (by induction) $(z,i)$ satisfies (\ref{hve}) for every $i \in I$. Therefore $h_i(z)=\min H$ for every $i$, hence $z \not\in X$, a contradiction.
Therefore $X=\emptyset$, hence $h_{i}(x)=h_{i+1}(x)$ for every $x \in A$ and $i \geq k$, thus $h_k=h_{k+1}$.
\end{proo}

\subsection{Valuation by measure}
Let us go back to the situation of significance domains on a space with a measure. Given $M \subset A$ and $x \in M'=A \setminus M$, let $\mu(x,M)=m(S_{x,M'})$ where $S_{x,M'}=\{y\in {\bf X}|(\forall a \in M') \: \delta(y,a)\geq \delta(y,x)\}$. Then one can easily see that the property (\ref{izoton}) is satisfied and the induced evolution of valuations is $(f_i)_{i \in \Nn_{0}}$. Then the Theorem \ref{stop} can be directly applied to this situation:
\begin{theo}
Let $X$ be a set with a metric $\delta$ and a measure $\mu$ and an initial valuation $h$ on $A$. The induced the evolution $(f_i)_{i\in \Nn_{0}}$ of valuation stops.
\end{theo}
Hence the system of valuations converges (in finite a number of steps) to a valuation. The limit valuation can be interpreted as a {\em potential} of a point w.r.t the initial valuations.

\end{document}